\documentclass[number,citesort,dvips]{arxstspdf}
\usepackage{graphicx}
\usepackage{flushend}
\usepackage{stfloats}

\volume{25}
\issue{3}
\pubyear{2010}
\firstpage{275}
\lastpage{288}
\doi{10.1214/10-STS335}

\makeatletter
\newcommand{\Ropt}{{R}_{\mathrm{{opt}}}}
\newcommand{\Rave}{{R}_{\mathrm{{ave}}}}
\newcommand{\Rmax}{{R}_{\mathrm{{max}}}}
\newcommand{\Rbeta}{{R}_{\beta\mbox{-}\mathrm{skel}}}
\newproclaim{definition}{Definition}
\newproclaim{remarks}{Remarks}
\newtheorem{Lemma}{Lemma}
\setattribute{ead}{sep}{; }
\makeatother

\begin{document}
\begin{frontmatter}

\title{Connected Spatial Networks over Random Points and a Route-Length
Statistic}
\runtitle{Connected Networks over Random Points}

\begin{aug}
\author[a]{\fnms{David J.} \snm{Aldous}\corref{}\ead[label=e1]{aldous@stat.berkeley.edu}
\ead[label=u1,url]{www.stat.berkeley.edu/users/aldous}}\and
\author[b]{\fnms{Julian} \snm{Shun}\ead[label=e2]{julianshun@gmail.com}}
\runauthor{D. J. Aldous and J. Shun}

\address[a]{David J. Aldous is Professor,
Department of Statistics,
University of California,
367 Evans Hall \# 3860,
Berkeley, California 94720,
USA
(\printead{e1,u1}).}
\address[b]{Julian Shun is Graduate Student,
Machine Learning Department,
Carnegie Mellon University,
5000 Forbes Avenue,
Pittsburgh, Pennsylvania 15213, USA (\printead{e2}).}

\end{aug}

%
\begin{abstract}
We review mathematically tractable models for connected networks on
random points in the plane, emphasizing the class of
\textit{proximity graphs} which deserves to be better known to applied
probabilists and statisticians.
We introduce and motivate a particular statistic $R$ measuring
shortness of routes in a network.
We illustrate, via Monte Carlo in part, the trade-off between
normalized network length and $R$ in a one-parameter family of
proximity graphs.
How close this family comes to the optimal trade-off over all possible
networks remains an intriguing open question.

The paper is a write-up of a talk developed by the first author during
2007--2009.

\end{abstract}
%
%

\end{frontmatter}
%

\section{Introduction}
\label{sec-INT}

The topic called \textit{random networks} or \textit{complex
networks} has
attracted huge
attention over the last 20 years.
Much of this work focuses on examples such as social networks or WWW
links, in which
edges are not closely constrained by two-dimensional geometry.
In contrast, in a \textit{spatial network} not only are vertices
and edges situated in two-dimensional space, but also it is actual
distances, rather than number of edges, that are of interest.
To be concrete, we visualize idealized inter-city road networks, and a
feature of
interest is the (minimum) route length between two given cities.
Because we work only in two dimensions, the word \textit{spatial} may be
misleading, but equally the word \textit{planar} would be misleading
because we do not require networks to be planar graphs (if edges cross,
then a junction is created).

Our major purpose is to draw the attention of readers from
the applied probability and statistics communities to a particular
class of spatial
network models.
Recall that the most studied network model, the \textit{random geometric
graph} \cite{penrose-RGG}
reviewed in Section \ref{sec-geom}, does not permit
both connectivity and bounded normalized length in the $n \to\infty$ limit.
An attractive alternative is the class of \textit{proximity graphs},
reviewed in
Section \ref{sec-proximity},
which in the deterministic case have been studied within computational geometry.
These graphs are always connected.
Proximity graphs on \textit{random} points have been studied in only a
few papers,
but are potentially interesting for many purposes other than the
specific ``short route lengths'' topic of this paper
(see Section \ref{sec-specu}).
One could also
imagine constructions which depend on points having specifically the
Poisson point process distribution,
and one novel such network, which we name the \textit{Hammersley network},
is described in Section \ref{sec-Hamm}.

Visualizing idealized road networks, it is natural to take total
network length
as the ``cost'' of a network, but what is the corresponding
``benefit''?
Primarily we are interested in having short route lengths.
Choosing an appropriate statistic to measure the latter turns
out to be rather subtle, and the (only) technical innovation of this
paper is
the introduction (Section \ref{sec-Rs}) and motivation of a specific
statistic $R$ for measuring the effectiveness of
a network in providing short routes.

In the theory of spatial networks over random points,
it is a challenge to quantify the trade-off between network length
[precisely, the normalized length $L$
defined at (\ref{theta-def})]
and route length efficiency statistics such as $R$.
Our particular statistic $R$ is not amenable to explicit calculation
even in comparatively tractable models,
but in Section \ref{sec-LET} we present the results from Monte Carlo
simulations.
In particular, Figure \ref{fig7} shows the trade-off for the particular \textit{
$\beta$-skeleton} family of proximity graphs.

Given a normalized network length $L$,
for any realization of cities
there is some network of normalized length $L$ which minimizes $R$.
As indicated in Section \ref{sec-limits},
by general abstract mathematical arguments, there must exist a
deterministic function
$\Ropt(L)$ giving
(in the
``number of cities $\to\infty$'' limit under the random model)
the minimum value of $R$ over all possible networks of normalized
length $L$.
An intriguing open question is as follows:
\[
\begin{tabular}{p{200pt}@{}}
how close are the values $\Rbeta(L)$ from the $\beta$-skeleton
proximity graphs to the optimum values
$\Ropt(L)$?
\end{tabular}
\]
As discussed in Section \ref{sec-curve}, at first sight it looks easy
to design heuristic algorithms for networks which should improve over
the $\beta$-skeletons, for example, by introducing Steiner points, but
in practice we have not succeeded in doing so.

This paper focuses on the random model for city positions because it
seems the natural setting for theoretical study.
As a complement,
in \cite{me-spatial-2} we give empirical data for the values of $(L,R)$
for certain real-world networks
(on the 20 largest cities, in each of 10 US States).
In \cite{me-spatial-3} we give analytic results and bounds on the
trade-off between $L$ and the mathematically more tractable
\textit{stretch} statistic $\Rmax$ at (\ref{def-Rmax}), in both worst-case
and random-case settings for city positions.
Let us also point out a (perhaps) nonobvious insight discussed in
Section \ref{sec-max23}:
in designing networks to be efficient in the sense of providing short routes,
the main difficulty is providing short routes between city-pairs at a
specific distance (2--3 standardized units) apart, rather than between
pairs at a large distance apart.

Finally, recall this is a nontechnical account.
Our purpose is to elaborate verbally the ideas outlined above; some
technical aspects will be pursued elsewhere.

\section{Models for Connected Spatial Networks}
\label{sec-connect}
There are several conceptually different ways of defining networks on
random points in the plane.
To be concrete, we call the points \textit{cities}; to be consistent about
language, we regard $x_i$ as
the \textit{position} of city $i$ and represent network edges as line
segments $(x_i,x_j)$.

First (Sections \ref{sec-geom}--\ref{sec-proximity}) are schemes which
use deterministic rules to define edges for an arbitrary deterministic
configuration of cities; then one just applies these rules to a random
configuration.
Second, one can have random rules for edges in a deterministic
configuration
(e.g., the probability of an edge between cities $i$ and $j$ is a
function of Euclidean distance $d(x_i,x_j)$, as in popular
\textit{small worlds} models \cite{newman-survey}), and again apply
to a
random configuration.
Third, and more subtly, one can have constructions that depend on the
randomness model for city positions---Section \ref{sec-Hamm} provides
a novel example.

We work throughout with reference to Euclidean distance $d(x,y)$ on the
plane, even though many models could be defined
with reference to other metrics (or even when the triangle inequality
does not hold, for the MST).

\subsection{The Geometric Graph}
\label{sec-geom}
In Sections \ref{sec-geom}--\ref{sec-proximity} we have an arbitrary
configuration
$\mathbf{x}= \{x_i\}$ of city positions, and a deterministic rule for
defining the edge-set $\mathcal{E}$.
Usually in graph theory one imagines a \textit{finite} configuration, but
note that everything makes sense for \textit{locally finite}
configurations too. Where helpful, we assume ``general position,'' so
that intercity distances $d(x_i,x_j)$ are all distinct.

\begin{figure*}[b]

\includegraphics{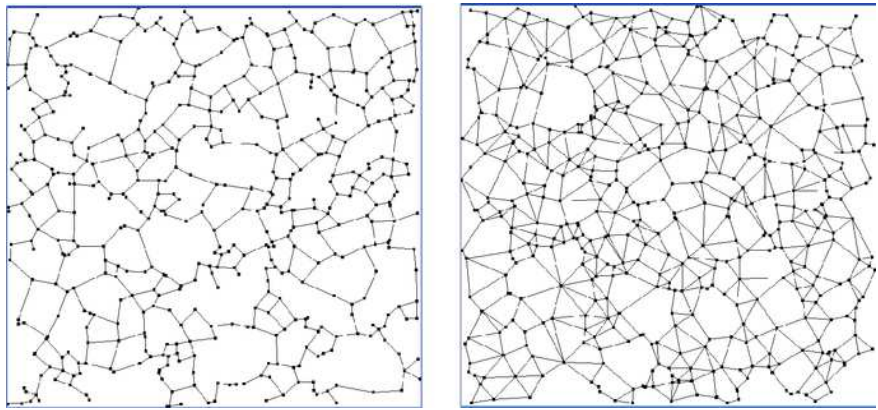}

  \caption{The relative neighborhood graph \textup{(left)} and Gabriel graph
\textup{(right)} on different realizations of $500$ random points.}\label{fig1}
\end{figure*}

For the \textit{geometric graph} one fixes $0 < c < \infty$ and defines
\[
(x_i,x_j) \in\mathcal{E}\quad \mbox{iff}\quad  d(x_i,x_j) \leq c .
\]
For the \textit{$K$-neighbor graph} one fixes $K \geq1$ and defines
\[
\begin{tabular}{p{200pt}@{}}
$(x_i,x_j) \in\mathcal{E}$ iff $x_i$ is one of the $K$ closest
neighbors of $x_j$, or $x_j$ is one of the $K$ closest neighbors of $x_i$.
\end{tabular}
\]
A moment's thought shows these graphs are in general not connected,
so we turn to models which are ``by construction'' connected.
We remark that the connectivity threshold $c_n$ in the \textit{finite}
$n$-vertex model of the random geometric graph has been studied
in detail---see Chapter 13 of \cite{penrose-RGG}.

\subsection{A Nested Sequence of Connected Graphs}
\label{sec-nested}
The material here and in the next section was developed in graph theory
with a view toward
algorithmic applications
in computational geometry and pattern recognition.
The 1992 survey \cite{jaromczyk} gives the history of the subject and
116 citations.
But everything we need is immediate from the (careful choice of) definitions.
On our arbitrary configuration $\mathbf{x}$ we can define four graphs whose
edge-sets are nested as follows:
%
\begin{equation}
\qquad \mbox{MST} \subseteq
\mbox{relative n'hood}
\subseteq
\mbox{Gabriel}
\subseteq
\mbox{Delaunay}.
\label{inclusion}
\end{equation}

Here are the definitions (for MST and Delaunay, it is easy to check
these are equivalent to more familiar definitions).
In each case, we write the criterion for an edge $(x_i,x_j)$ to be present:

\begin{itemize}
\item
\textit{Minimum spanning tree (MST) \textup{\cite{graham-hell}}.}
There does not exist a sequence $i = k_0,k_1,\ldots, k_m = j$ of cities
such that
\begin{eqnarray*}
&&\max(d(x_{k_0},x_{k_1}), d(x_{k_1},x_{k_2}), \ldots,
d(x_{k_{m-1}},x_{k_m})) \\
&&\quad< d(x_i,x_j).
\end{eqnarray*}
\item
\textit{Relative neighborhood graph.}
There does not exist a city $k$ such that
\[
\max(d(x_i,x_k), d(x_k,x_j)) < d(x_i,x_j).
\]
\item
\textit{Gabriel graph.}
There does not exist a city inside the disc whose diameter is the line
segment from $x_i$ to~$x_j$.
\item
\textit{Delaunay triangulation \textup{\cite{MR1686530}}.}
There exists some disc, with $x_i$ and $x_j$ on its boundary, so that
no city is inside the disc.
\end{itemize}
The inclusions (\ref{inclusion}) are immediate from these definitions.
Because the MST (for a finite configuration) is connected, all these
graphs are connected.

Figure \ref{fig1} illustrates the relative neighborhood and Gabriel graphs.
Figures for the MST and the
Delaunay triangulation can be found online at
\href{http://www.spss.com/research/wilkinson/Applets/edges.html}{http://www.}\break
\href{http://www.spss.com/research/wilkinson/Applets/edges.html}{spss.com/research/wilkinson/Applets/edges.html}.

Constructions such as the relative neighborhood and Gabriel graphs have
become known loosely as
\textit{proximity graphs} in \cite{jaromczyk} and subsequent literature,
and we next take the opportunity to turn an implicit definition in the
literature into an explicit definition.

\subsection{Proximity Graphs}
\label{sec-proximity}
Write $v_{-}$ and $v_+$ for the points
$(-{ \frac{1}{2}},0)$ and $({ \frac{1}{2}},0)$.
The \textit{lune} is the intersection of the open discs of radii~$1$
centered at $v_{-}$ and $v_+$.
So $v_{-}$ and $v_+$ are not in the lune but are on its boundary.
Define a \textit{template} $A$ to be a subset of ${\mathbb{R}}^2$
such that:
\begin{longlist}[(iii)]
\item[(i)] $A$ is a subset of the lune.
\item[(ii)] $A$ contains the open line segment $(v_{-},v_+)$.
\item[(iii)] $A$ is invariant under the ``reflection in the $y$-axis'' map
$\operatorname{Reflect}_x (x_1,x_2) = (-x_1,x_2)$
and the ``reflection in the $x$-axis'' map $\operatorname{Reflect}_y
(x_1,x_2) = (x_1,\break  -x_2)$.
\item[(iv)] $A$ is open.
\end{longlist}

For arbitrary points $x, y$ in ${\mathbb{R}}^2$, define $A(x,y)$ to be the
image of $A$ under the natural
transformation (translation, rotation and scaling) that takes
$(v_-,v_+)$ to $(x,y)$.

\begin{definition*}
Given a template $A$ and a locally finite set $\mathcal{V}$ of
vertices, the
associated \textit{proximity graph} $G$ has edges defined by,
for each $x, y \in\mathcal{V}$,
\[
\begin{tabular}{p{200pt}@{}}
$(x,y)$ is an edge of $G$ iff $A(x,y)$ contains no vertex of
$\mathcal{V}$.
\end{tabular}
\]
From the definitions:
\begin{itemize}
\item if $A$ is the lune, then $G$ is the relative neighborhood graph;
\item if $A$ is the disc centered at the origin with radius $1/2$, then
$G$ is the Gabriel graph.
\end{itemize}
But the MST and Delaunay triangulation are not instances of
proximity graphs.
\end{definition*}

Note that replacing $A$ by a subset $A^\prime$ can only introduce
extra edges.
It follows from (\ref{inclusion}) that 
the proximity graph is always connected.
The Gabriel graph is planar.
But if $A$ is not a superset of the disc centered at the origin with
radius $1/2$, then $G$ might not be a subgraph of the Delaunay
triangulation, and in this case edges may cross, so $G$ is not planar
(e.g., if the vertex-set is the four corners of a square, then the
diagonals would be edges).

For a given configuration $\mathbf{x}$, there is a collection of proximity
graphs indexed by the template $A$,
so by choosing a monotone one-parameter family of templates, one gets a
monotone one-parameter family
of graphs,
analogous to
the one-parameter family $\mathcal{G}_c$ of geometric graphs.
Here is a popular choice \cite{kirk-radke} in which $\beta= 1$ gives
the Gabriel graph and $\beta= 2$ gives the
relative neighborhood graph.

\begin{definition*}[(The \textit{$\beta$-skeleton} family)]
(i) For $0<\beta< 1$ let $A_\beta$ be the intersection of the two open
discs of radius
$(2 \beta)^{-1}$ passing through $v_{-}$ and $v_{+}$.

(ii) For $1 \leq\beta\leq2$ let $A_\beta$ be the intersection of the two
open discs of radius
$\beta/2$ centered at $(\pm(\beta-1)/2,0)$.
\end{definition*}

\subsection{Networks Based on Powers of Edge-Lengths}
\label{sec-kumar}
It is not hard to think of other ways to define one-parameter families
of networks.
Here is one scheme used in, for example, \cite{kumar02}.
Fix $1 \leq p < \infty$.
Given a configuration $\mathbf{x}$, and a route (sequence of vertices)
$x_0,x_1,\ldots,x_k$, say, the cost of the route is the sum of $p$th
powers of the step lengths.
Now say that a pair $(x,y)$ is an edge of the network $\mathcal{G}_p$
if the
cheapest route from $x$ to $y$ is the one-step route.
As $p$ increases from~$1$ to $\infty,$ these networks decrease from the
complete graph to the MST.
Moreover, for $p \geq2$ the network $\mathcal{G}_p$ is a subgraph of the
Gabriel graph.

\subsection{The Hammersley Network}
\label{sec-Hamm}
There is a quite
separate recent literature
in theoretical probability
\cite{Poisson-Matching,MR1961286}
defining structures such as trees and matchings directly on the
infinite Poisson point process.
In this spirit, we observe that the
\textit{Hammersley process}
studied in \cite{me71}
can be used
to define a new network on the infinite Poisson point process,
which we name the \textit{Hammersley network}.
This network is designed to have the feature that
each vertex has exactly $4$~edges,
in directions NE (between North and East), NW, SE and SW.
The conceptual difference from the networks in the previous section is
that there is not such a simple ``local'' criterion for whether a
potential edge $(x_i,x_j)$ is in the network.
And edges cross, creating junctions.

For a picturesque description, imagine one-eyed frogs sitting on an
infinitely long, thin log, each being able to see only the part of the
log to their left before the next frog.
At random times and positions (precisely, as a space--time Poisson point
process of rate $1$) a fly lands on the log, at which instant the
(unique) frog which can see it jumps left to the fly's position and
eats it.
This defines a continuous time Markov process (the Hammersley process)
whose states are the configurations of positions of all the frogs.
There is a \textit{stationary} version of the process in which, at each
time, the positions of the frogs form a Poisson (rate $1$) point
process on the line.

\begin{figure}[b]

\includegraphics{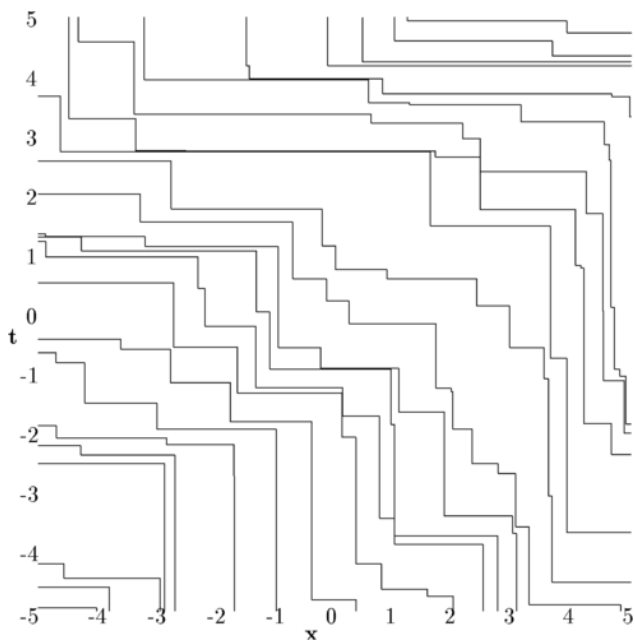}

  \caption{Space--time trajectories in Hammersley's process.}\label{fig2}
\end{figure}

Now consider the space--time trajectories of all the frogs, drawn with
time increasing upward on the page.
See Figure \ref{fig2}.
For each frog, the part of the trajectory between the completions of
two successive jumps consists of an upward edge
(the frog remains in place as time increases) followed by a leftward
edge (the frog jumps left).

\begin{figure}

\includegraphics{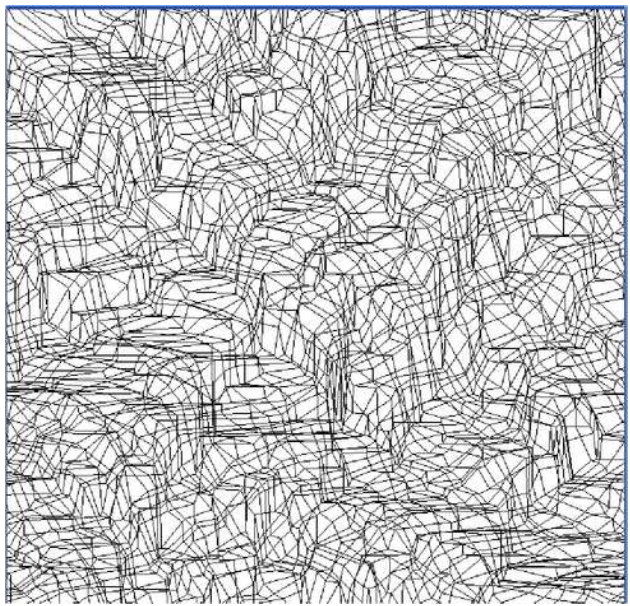}

  \caption{The Hammersley network on 2500 random points.}\label{fig3}
\end{figure}

Reinterpreting the time axis as a second space axis, and introducing
compass directions,
that part of the trajectory becomes a North edge followed by a West edge.
Now replace these two edges by a single North-West straight edge.
Doing this procedure for each frog and each pair of successive jumps,
we obtain
a collection of NW paths, that is, a network in which each city (the
reinterpreted space--time random points) has an edge to the NW and an
edge to the SE.
Finally, we repeat the construction with the same realization of the
space--time Poisson point process but with frogs jumping rightward
instead of leftward.
This yields a network on the infinite Poisson point process, which we
name the \textit{Hammersley network}. See Figure \ref{fig3}.

\begin{remarks*}
(a) To draw the Hammersley network on random points in a \textit{finite
square}, one needs
external randomization to give the initial (time $0$) frog positions,
in fact, two independent randomizations for the leftward and the
rightward processes.
So to be pedantic, one gets a \textit{random} network over the given
realization of cities.
However, one can deduce from the theoretical results in \cite{me71}
that the external randomization has effect only near the boundary of
the square.\vspace*{-6pt}
\begin{longlist}[(b)]
\item[(b)] The property that
each vertex has exactly\break $4$~edges,
in directions NE (between North and East), NW, SE and SW,
is immediate from the construction.
Note, however, that while adjacent NW space--time trajectories in Figure
\ref{fig2} do not cross, the corresponding diagonal roads
in the Hammersley network may cross, so it is not a planar graph,
though this has only negligible effect on route lengths.

\item[(c)] Intuition, confirmed by Figure \ref{fig7} later, says that the
Hammersley network is not very efficient as a road network.
It serves to demonstrate that there do exist random networks other
than the familiar ones, and
provides an instance where
imposing deterministic constraints (the four edges, in this case) on a
random network
makes it much less efficient. How general a phenomenon is this?
\end{longlist}
\end{remarks*}

\subsection{Normalized Length}
The notion of
\textit{normalized network length}
$L$
is most easily visualized in the setting of an infinite deterministic network
which is ``regular'' in the sense of consisting of a repeated pattern.
First choose the unit of length so that cities have an average density
of one per unit area.
Then define
%
%
\begin{eqnarray}\label{theta-def}
\qquad L &=& \mbox{average network length per unit area,}\\
\bar{\Delta} &=& \mbox{average degree (number of incident edges)}\label{bard-def}
\nonumber
\\[-8pt]
\\[-8pt]
\nonumber
&&{} \mbox{of cities.}
\end{eqnarray}

\begin{figure*}

\includegraphics{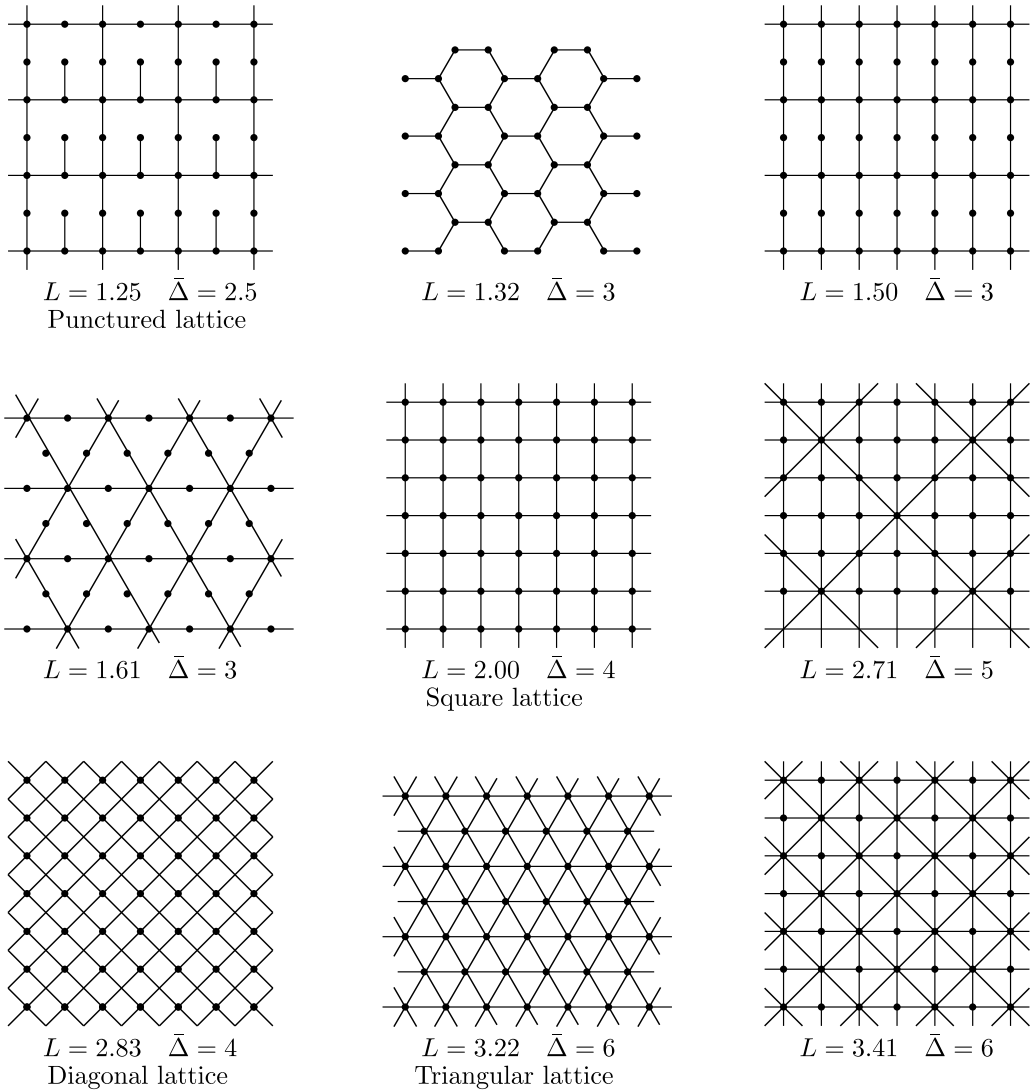}

  \caption{Variant square, triangular and hexagonal lattices.
Drawn so that the density of cities is the same in each
diagram, and ordered by value of $L$.}\label{fig4}
\end{figure*}

Figure \ref{fig4} shows the values of $L$ and $\bar{\Delta}$ for some simple
``repeated pattern'' networks.
Though not directly relevant to our study of the random model, we find
Figure \ref{fig4} helpful for two reasons: as intuition for the interpretation
of the different numerical values of $L$, and because we can make very
loose analogies (Section \ref{sec-analogs}) between particular networks
on random points and particular deterministic networks.

\section{Normalized Length and route-Length Efficiency}
\subsection{The Random Model}
For the remainder of the paper we work with ``the random model'' for
city positions.
The \textit{finite model} assumes $n$ random vertices (cities) distributed
independently and uniformly in a square of area $n$.
The \textit{infinite model} assumes the Poisson point process of rate $1$
(per unit area) in the plane.
The quantities $L, \bar{\Delta}$ above and $R$ below that we discuss
may be interpreted as exact values in the infinite model or as $n \to
\infty$ limits in the finite model; see Section \ref{sec-limits}.
We use the word \textit{normalized} as a reminder of the ``density $1$''
convention---we choose the normalized unit of distance to make cities
have average density $1$ per unit area.
After this normalization, $L$ is the average network length per unit area.

\subsection{The Route-Length Efficiency Statistic $R$}
\label{sec-Rs}
In designing a network, it is natural to regard total length as a ``cost''.
The corresponding ``benefit'' is having short routes between cities.
Write $\ell(i,j)$ for the route length (length of shortest path)
between cities
$i$ and $j$ in a given network, and $d(i,j)$
for Euclidean distance between the cities.
So $\ell(i,j) \geq d(i,j)$, and we write
\[
r(i,j) = { \frac{\ell(i,j)}{d(i,j)}} - 1
\]
so that ``$r(i,j) = 0.2$'' means that route length is $20\%$ longer than
straight line distance.
With $n$ cities we get ${n \choose2}$ such numbers $r(i,j)$; what is a
reasonable way to combine these into a single statistic? Two natural
possibilities are as follows:
%
%
\begin{eqnarray}\label{def-Rmax}
\Rmax&:=&\max_{j \neq i} r(i,j),
\nonumber
\\[-8pt]
\\[-8pt]
\nonumber
\Rave&:=& \operatorname{ave}_{(i,j)} r(i,j),
\end{eqnarray}
where $\operatorname{ave}_{(i,j)}$ denotes average over all distinct pairs $(i,j)$.
The statistic $\Rmax$ has been studied in the context of the design of
geometric spanner networks \cite{MR2289615} where it is called the
\textit{stretch}.
However, being an ``extremal'' statistic $\Rmax$ seems unsatisfactory
as a descriptor of real world networks---for instance, it seems
unreasonable to characterize the UK rail network as inefficient
simply because there is no very direct route between Oxford and Cambridge.

The statistic $\Rave$ has a more subtle drawback.
Consider a network consisting of:
\begin{itemize}
\item the minimum-length connected network (Steiner tree) on given cities;
\item
and a superimposed sparse collection of randomly oriented lines
(a \textit{Poisson line process} \cite{MR895588}).
\end{itemize}

\begin{figure}

\includegraphics{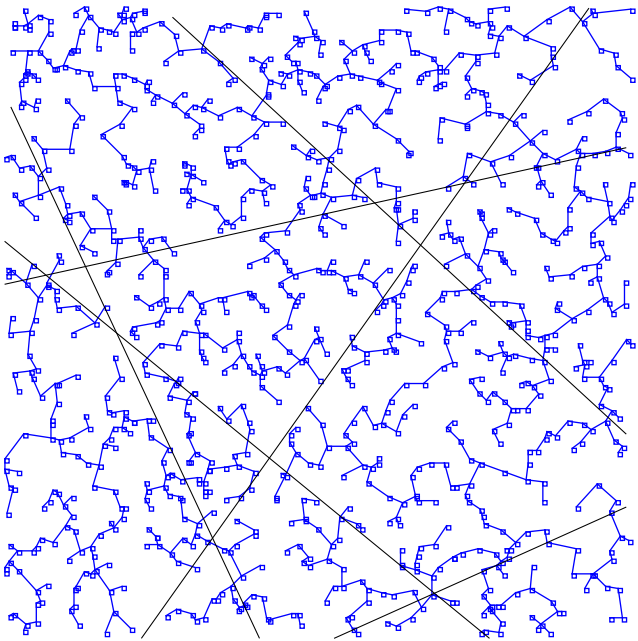}

  \caption{Efficient or inefficient?
$\Rave$ would judge this network efficient in the $n \to
\infty$ limit.}\label{fig5}
\end{figure}

\noindent
See Figure \ref{fig5}.
By choosing the density of lines to be sufficiently low, one can make the
normalized network length be arbitrarily close to the minimum needed
for connectivity.
But it is easy to show (see \cite{me116} for careful analysis and a
stronger result) that one can construct such
networks so that $\Rave\to0$ as $n \to\infty$.
Of course no one would build a road network looking like Figure~\ref{fig5} to
link cities, because there are many pairs of nearby cities with only
very indirect routes between them.
The disadvantage of $\Rave$ as a descriptive statistic is that (for
large $n$) most city-pairs are far apart, so the fact that a given
network has a small value of $\Rave$ says nothing about
route lengths between \textit{nearby} cities.

We propose a statistic $R$ which is intermediate between
$\Rave$ and $\Rmax$.
First consider (see discussion below for details)
\begin{eqnarray*}
\rho(d) &:=&
\mbox{mean value of $r(i,j)$ over}\\
&&{}\mbox{city-pairs with $d(i,j) = d$}
\end{eqnarray*}
and then define
%
%
\begin{equation}
R := \max_{0 \leq d < \infty} \rho(d) . \label{R-def}
\end{equation}
In words, $R = 0.2$ means that on every scale of distance, route
lengths are on average at
most $20\%$ longer than straight line distance.

On an intuitive level, $R$ provides a
sensible and interpretable way to compare efficiency of different
networks in providing
short routes.
On a technical level,
we see two advantages and one disadvantage
of using $R$ instead of
$\Rave$.

\textit{Advantage} 1.
Using $R$ to measure efficiency, there is a meaningful $n \to\infty$
limit for the network length/\break efficiency trade-off [the function $\Ropt
(L)$ discussed in Section \ref{sec-limits}], and so, in particular, it
makes sense to
compare the values of $R$ for networks with different~$n$.

\textit{Advantage} 2.
A more realistic model for traffic would posit that volume of traffic
between two cities varies as a power-law $d^{-\gamma}$ of distance $d$,
so that in calculating $\Rave$ it would be more realistic to weight by
$d^{-\gamma}$. This means that the optimal network, when using $\Rave$
as optimality criterion, would depend on $\gamma$. Use of $R$ finesses
this issue; the value of $\gamma$ does not affect~$R$.
A related issue is that volume of traffic between two cities should
depend on their populations.
Intuitively, incorporating random population sizes should make the
optimal $R$
\textit{smaller} because the network designer can create shorter routes
between larger cities.
We see this effect in data \cite{me-spatial-2}; $R$ calculated via
population-weighting is typically slightly smaller.
But we have not tried theoretical study.

\textit{Disadvantage}.
The statistic $R$ is tailored to the infinite model, in which it makes
sense to consider
two cities at exactly distance $d$ apart
(then the other city positions form a Poisson point process).
For finite $n$ we need to discretize.
For the empirical data in
\cite{me-spatial-2}, where $n = 20$,
we average over intervals of width $1$ unit
(recall the unit of distance is taken such that the density of cities
is $1$
per unit area),
that is, for $ d = 1,2, \ldots, 5,$ we calculate
\begin{eqnarray}\label{tilde-R}
\tilde{\rho}(d)&:=&
\mbox{mean value of $r(i,j)$ over city-pairs} \nonumber\\
&&{}\mbox{with }d - { \tfrac
{1}{2}} <
d(i,j) < d + \tfrac{1}{2},
\\
\tilde{R} &:=& \max_{1 \leq d < \infty} \tilde{\rho}(d)\nonumber
\end{eqnarray}
and use $\tilde{R}$ as proxy for $R$.
For larger $n$ we can use shorter intervals.
Thus, there is, in principle,
a certain fuzziness to the notion of $R$ for finite networks, and,
in particular,
it is not clear how to assign a value of $R$ to regular networks such
as those
in Figure \ref{fig4}.
But in practice, for networks we have studied on real-world data and on
random points, this is not a problem, as explained next.

\begin{figure*}

\includegraphics{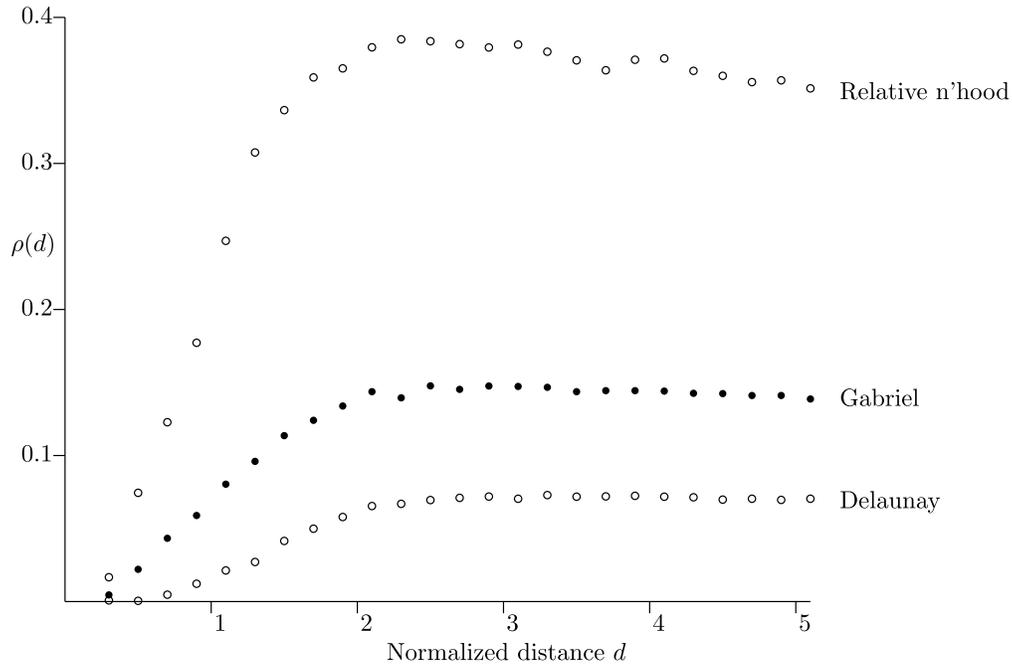}

  \caption{The function $\rho(d)$ for three theoretical networks on
random cities.
Irregularities are Monte Carlo random variation.}\label{fig6}
\end{figure*}

\subsection{Characteristic Shape of the Function $\rho(d)$}
\label{sec-max23}
For the connected networks on random points (excluding the Hammersley
network) we are discussing,
the function $\rho(d)$ has a characteristic shape
(see Figure~\ref{fig6})
attaining its maximum between $2$ and $3$ and slowly decreasing thereafter.
We suspect that ``this characteristic shape holds for any reasonable
model,'' but we do not know how to turn that phrase into a precise conjecture.
Note that ``smoothness near the maximum'' implies
that any calculated
value $\tilde{R}$ at (\ref{tilde-R}) is quite insensitive to the choice
of discretization.

This characteristic shape has a common-sense interpretation.
Any efficient network will tend to place roads directly between
unusually close city-pairs, implying that $\rho(d)$ should be small for
$d < 1$.
For large $d$
the presence of multiple alternate routes helps prevent $\rho(d)$ from growing.
At distance $2 - 3$ from a typical city $i$ there will be about $\pi
3^2 - \pi2^2 \approx16$ other cities~$j$.
For some of these $j$ there will be cities $k$ near the straight line
from $i$ to $j$, so the network designer can create roads from $i$ to
$k$ to $j$.
The difficulty arises where there is no such intermediate city $k$:
including a direct road $(x_i,x_j)$ will increase $L$, but not
including it will increase $\rho(d)$ for
$2<d<3$.

Thus, Figure \ref{fig6} offers a minor insight into spatial network design:
that it is city pairs at normalized distance $2 - 3$ specifically that
enforce the constraints on efficient network design.

The characteristic shape---at least, the flatness over $2 \le d \le
5$---is also visible in the real-world data \cite{me-spatial-2}.

For the Hammersley network, the graph of $\rho(d)$ is quite different;
$\rho(d)$ increases to a maximum of $0.35$ around $d = 0.8$ and then decreases
more steeply to a value of $0.21$ at $d = 5$.
This arises from the particular structure (from each city there is one
road in each quadrant) resembling the deterministic
``diagonal lattice'' of Figure \ref{fig4}, in which the route between some nearby
pairs will be via two diagonal roads and a junction.

\section{Length-Efficiency Trade-Off for Tractable Networks}
\label{sec-LET}
Recall that our overall theme
is the trade-off between network length and route-length efficiency,
and that in this paper we focus on $n \to\infty$ limits in the random model
and the particular statistics $L$ and $R$.

\begin{figure*}

\includegraphics{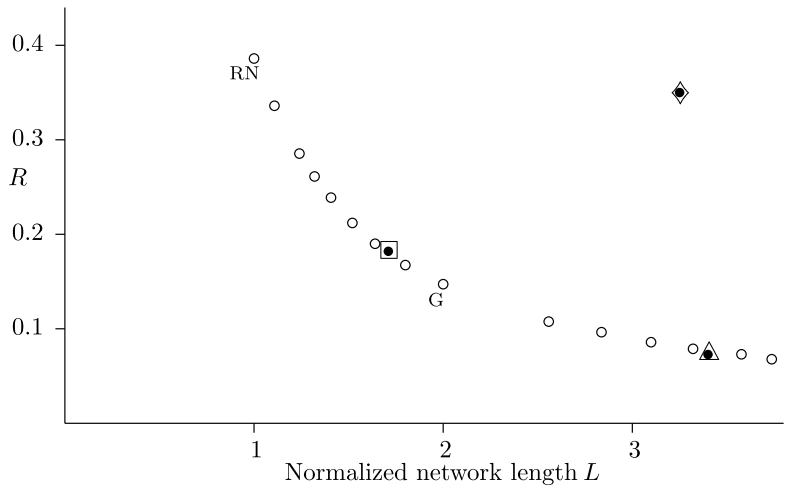}

  \caption{The normalized network length $L$ and the route-length efficiency
statistic $R$ for certain
networks on random points. The $\circ$ show the beta-skeleton family,
with RN the relative neighborhood graph and G the Gabriel graph. The
$\bullet$ are special models: $\triangle$ shows the
Delaunay triangulation, $\Box$ shows the network $\mathcal{G}_2$ from Section
\protect\ref{sec-kumar} and $\diamondsuit$ shows the Hammersley network.}\label{fig7}
\end{figure*}

The models described in Section \ref{sec-connect}
are ``tractable'' in the specific sense that one can find exact analytic
formulas for normalized length $L$.
Unfortunately $R$ is not amenable to analytic calculation, 
and we resort to Monte Carlo simulation to obtain values for $R$.
Table \ref{tab1} and Figure \ref{fig7}
show the values of $(L,R)$ in the models.
We explain below how the values of $L$ are calculated.


\begin{table}[b]
\caption{Statistics of tractable networks on random
points}\label{tab1}
\begin{tabular*}{\columnwidth}{@{\extracolsep{\fill}}l@{\quad}ccc@{}}
\hline
\textbf{Network} & $\bolds{L}$ & $\bolds{\bar{\Delta}}$ & $\bolds{R}$\\
\hline
Minimum spanning tree& 0.633 & 2\phantom{00.} & $\infty$
\\
Relative n'hood & 1.02\phantom{0} & $2.56$ & 0.38
\\
Gabriel & 2\phantom{00.0} & $4$\phantom{00.} & 0.15
\\
Hammersley & 3.25\phantom{0} & $4$\phantom{00.} & 0.35
\\
Delaunay & 3.40\phantom{0} & 6\phantom{00.} & 0.07\\
\hline
\end{tabular*}
\tabnotetext[]{}{\textit{Notes}:
Integer values are exact.
Recall $L$ is normalized length (\protect\ref{theta-def}), $\bar{\Delta}$ is
average degree (\protect\ref{bard-def}) and $R$ is our route-length statistic
(\protect\ref{R-def}).}
\end{table}

\textit{Notes on Table} \ref{tab1}.
(a) Values of $R$ from our simulations with $n = 2500$.

(b) Value of $L$ for MST from Monte Carlo \cite{cortina}.
In principle, one can calculate arbitrarily close bounds \cite{AB92},
but apparently this has never been carried out.
Of course, $\bar{\Delta} = 2$ for any tree.

(c) The Gabriel graph and the relative neighborhood graph
fit the assumptions of Lemma \ref{Lscale} with
$c = \pi/4$ and $c = \frac{2\pi}{3} - \frac{\sqrt{3}}{4}$, respectively,
and their table entries for $L$ and $\bar{\Delta}$ are obtained from
Lemma \ref{Lscale},
as are the values for $\beta$-skeletons in Figure \ref{fig7}.

(d) For the Hammersley network, every degree equals $4$, so $L = 2\times{}$(mean edge-length).
It follows from theory \cite{me71} that a typical edge, say, NE from $(x,y)$,
goes to a city at position $(x+ \xi_x, y + \xi_y)$, where $\xi_x$ and
$\xi_y$ are independent with $\operatorname{Exponential}(1$) distribution. So mean
edge-length equals
%
%
\begin{equation}
\int_0^\infty\!\!\!\int_0^\infty\sqrt{x^2 + y^2} e^{-x-y} \, dx\,dy
\approx1.62.
\label{Hamm-integral}
\end{equation}

(e) For any triangulation, $\bar{\Delta} = 6$ in the infinite model.
For the Delaunay triangulation, $L = ES$ where $S$ is the perimeter
length of a typical cell,
and it is known (\cite{miles70}, page 113) that $ES = \frac{32}{3\pi}$.
Note \cite{MR0660693} that the Delaunay triangulation is in general
\textit{not} the minimum-length triangulation.
Our simulation results in Figure~\ref{fig6} for $\rho(d)$ for the Delaunay
triangulation are roughly consistent with a simulation result in
\cite{BTZ00} saying that
$\rho(65) \approx0.05$.

\subsection{A Simple Calculation for Proximity Graphs}
Let us give an example of an elementary calculation for proximity
graphs over random points.
\begin{Lemma}
\label{Lscale}
For a proximity graph with template $A$ on the Poisson point process,
%
%
\begin{eqnarray}
L &=& { \frac{\pi^{3/2}}{4 c^{3/2}}},
\label{theta-scale}\\
\bar{\Delta} &=& { \frac{\pi}{c}},
\label{degree-scale}
\end{eqnarray}
where $c = \operatorname{area}(A)$.
\end{Lemma}

\begin{pf}
Take a typical city at position $x_0$.
For a city~$x$ at distance $s$ the chance that $(x_0,x)$ is an edge equals
$\exp(-cs^2)$ and so
\begin{eqnarray*}
\mbox{mean-degree} &=&
\int_0^\infty\exp(-cs^2)  2\pi s \, ds, \\
L &=&
{ \frac{1}{2}}
\int_0^\infty s \exp(-cs^2) 2\pi s \, ds .
\end{eqnarray*}
Evaluating the integrals gives (\ref{theta-scale}) and (\ref{degree-scale}).
\end{pf}

One can derive similar integral formulas for
other ``local'' characteristics, for example, mean density of triangles
and moments of vertex degree.
See \cite{MR937563,MR2257270,DGM09,matula} for a variety of such
generalizations and specializations.

\subsection{Other Tractable Networks}
\begin{figure}

\includegraphics{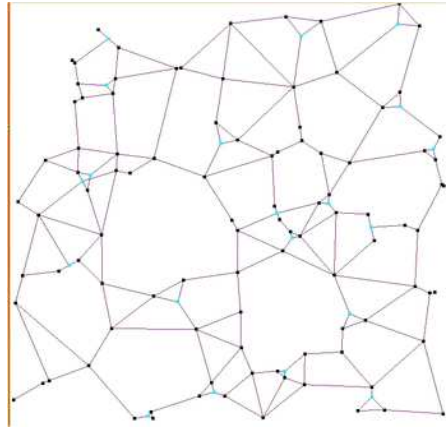}

  \caption{An ad hoc modification of the relative neighborhood graph, introducing
junctions.}\label{fig8}
\end{figure}

We do not know any
other ways of defining networks on random points which are both
``natural'' and are tractable in the sense
that one can find exact analytic formulas for $L$.
In particular, we know no tractable way of defining networks with
deliberate junctions as in Figure \ref{fig8}.
Note also that, while it is easy to make ad hoc modifications to the
geometric graph to ensure connectivity, these destroy tractability.
On the other hand, one can construct ``unnatural'' networks (see, e.g.,
\cite{me-spatial-3}) designed to permit calculation of $L$.

\section{Optimal Networks and $n \to\infty$ Limits}
\label{sec-limits}
\subsection{Tractable Models}
As mentioned earlier,
the quantities $L, \bar{\Delta}, R$ we discuss may be interpreted as
exact values in the infinite model or as $n \to\infty$ limits in the
finite model. To elaborate briefly, in a realization of the finite model
($n$ cities distributed independently and uniformly in a square of
area $n$), a network in Table \ref{tab1} has a normalized length $L_n = n^{-1}
\times\mbox{ (network length)}$ and an average degree $\bar{\Delta}_n$
which are random variables, but there is convergence (in probability
and in expectation)
%
%
\begin{equation}
L_n \to L,\quad  \bar{\Delta}_n \to\bar{\Delta}\quad  \mbox{as } n \to
\infty
\label{lims}
\end{equation}
to limit constants definable in terms of the analogous network on the
infinite model
(rate $1$ Poisson point process on the infinite plane).
For the proximity graphs or Delaunay triangulation, the network
definition applies directly to the infinite model and proof of (\ref
{lims}) is straightforward.
For the Hammersley network, (\ref{lims}) is implicit in \cite{me71},
and for the MST detailed arguments can be found in \cite{me57,MR1952000}.

\subsection{Optimal Networks}
We now turn to consideration of \textit{optimal} networks.
Given a configuration $\mathbf{x}$ of $n$ cities in the area-$n$ square,
and a value of $L$ which is greater than $n^{-1} \times\mbox{(length
of Steiner tree)}$,
one can define a number
%
%
\begin{eqnarray}\label{optimal}
\quad R_n(\mathbf{x}, L)& =&
\mbox{min of }\tilde{R} \mbox{ over all networks}
\nonumber
\\[-8pt]
\\[-8pt]
\nonumber
&&\mbox{on }\mathbf{x}
\mbox{ with normalized
length $\leq L$},
\end{eqnarray}
where $\tilde{R}$ is the discretized version (\ref{tilde-R}) calculated
using intervals of some suitable length $\delta_n$.
Applying this to a random configuration $\mathbf{X}$ in the finite model
gives, for each $L$,
a random variable
\[
\Xi_n(L) := R_n(\mathbf{X},L) .
\]
One intuitively expects convergence to some deterministic limit
%
%
\begin{equation}
\Xi_n(L) \to\Ropt(L) \quad \mbox{say, as } n \to\infty.
\label{Xi-con}
\end{equation}
The analogous result for $\Rmax$ will be proved carefully in \cite
{me-spatial-3}, and the same ``superadditivity''
argument could be used to prove (\ref{Xi-con}).
See \cite{MR1952000,steele97,yukich-book} for general background to
such results.
The point is that we do not have any explicit description of the optimal
[i.e., attaining the minimum in (\ref{optimal})]
networks in the finite or infinite models, so it seems very
challenging to prove the natural stronger supposition
that the finite optimal networks themselves converge (in some
appropriate sense) to a \textit{unique} infinite optimal network for which
the value $R = \Ropt(L)$ is attained.

\subsection{The Curve $\Ropt(L)$}
\label{sec-curve}
\begin{figure*}

\includegraphics{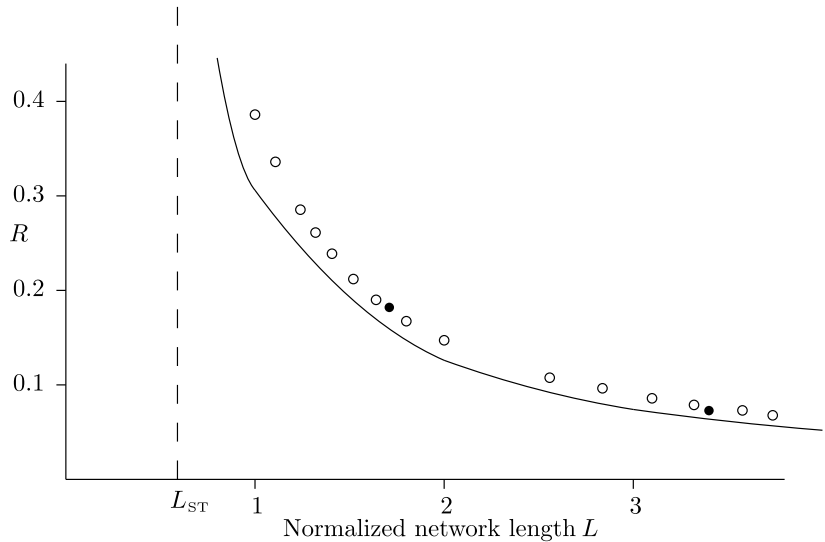}

  \caption{Speculative shape for the curve $\Ropt(L)$, with
$\circ$ and $\bullet$ values from tractable networks in Figure \protect\ref{fig7}.}\label{fig9}
\end{figure*}

Every possible network on the infinite Poisson point process defines a pair
$(L,R)$, and the curve $R = \Ropt(L)$ can be defied equivalently as the
lower boundary of the set of possible values
of $(L,R)$.
There is no reason to believe that proximity graphs are exactly
optimal, and, indeed, Figure \ref{fig7} shows that the Delaunay triangulation is
slightly more efficient than the corresponding $\beta$-skeleton.
But our attempts to do better by ad hoc constructions
(e.g., by introducing degree-$3$ junctions---see Figure \ref{fig8} for an example)
have been unsuccessful.
And, indeed, the fact that the
two special models in Figure \ref{fig7} lie close to the $\beta$-skeleton curve
lends credence to the idea that this curve is almost optimal.
We therefore speculate that the function $\Ropt$ looks something like
the curve in Figure \ref{fig9}, which we now discuss.

What can we say about $\Ropt(L)$? It is \textit{a priori} nonincreasing.
It is known \cite{yukich-book} that there exists a
\textit{Euclidean Steiner tree constant} $L_{\mathrm{ST}}$
representing the limit
normalized Steiner tree length in the random model, and clearly
$\Ropt(L) = \infty$ for $L < L_{\mathrm{ST}}$.
The facts
%
%
\begin{eqnarray}\label{decreases}
\Ropt(L) &<& \infty\quad \mbox{for all } L > L_{\mathrm{ST}};
\nonumber
\\[-8pt]
\\[-8pt]
\nonumber
\Ropt(L) &\to&0 \quad \mbox{as } L \to\infty
\end{eqnarray}
are not trivial to prove rigorously, but follow from the corresponding
facts for
$\Rmax$ proved in \cite{me-spatial-3}.
But we are unable to prove rigorously that $\Ropt(L)$ is \textit{strictly}
decreasing or that it is continuous.


\section{Final Remarks}
\label{sec-rem}

\subsection{Toy Models for Road Networks}
The idea of using proximity graphs as toy models for road networks has
previously been noted
\cite{kirk-radke} but not investigated very thoroughly.
It is an intuitively natural idea to a network designer:
whether or not to place a direct road from city $i$ to a nearby city
$j$ depends (partly) on whether some other city $k$ is close to the
line between them.

As observed by a referee, for the kind of models studied in this paper
we expect
route length $\ell(i,j)$ between distant cities to be roughly
proportional to graph distance (number of edges), which is a more
relevant quantity in some contexts.
However, when one considers design of optimal networks, replacing or
partially replacing
route length by graph distance leads to quite different optimal
networks \cite{me118,gastner-newman}.
For some other cost/benefit functionals leading to yet different
optimal networks see \cite{me112,barth06}.

\subsection{Rigorous Proof of Finite $R$ in Random Proximity Graphs}
\label{sec-rigor}
Table \ref{tab1} presented the Monte Carlo numerical value $\approx$0.38 of
$R$ for the relative neighborhood graph on random points.
From a rigorous viewpoint, the assertion that a random network has $R <
\infty$ is
essentially the assertion that
$\rho(d) = O(d)$ as $d \to\infty$. This is often nontrivial to prove.
A general sufficient condition for this property, which applies to the
relative neighborhood graph
(and hence all proximity graphs), is proved in \cite{me-spatial-4}.
The related fact that the limit $\lim_{d \to\infty} \rho(d)/d$ exists
is proved in
\cite{me-spatial-5}.

\subsection{Real-World Trade-Off Between Network Length and
Route-Length Efficiency}
Recall that our
central theme is seeking to quantify the trade-off between normalized
network length~$l$
and route-length efficiency $R$.
Figure \ref{fig9} suggests that for optimal networks the ``law of diminishing
returns'' sets in around $L = 2$
(for comparison, this is the value of $L$ corresponding to the square
grid network), in that
$\Ropt(L)$ decreases rapidly to around $0.13$ as $L$ increases to $2$
but decreases only slowly as $L$ increases further.
This suggests a kind of ``economic prediction'' for the lengths of
real-world networks which are perceived by users to be efficient
in providing short routes:
\[
\begin{tabular}{p{200pt}@{}}
the length of an efficient network linking $n$ cities in a region of
area $A$ will be roughly
$2 \sqrt{An}$.
\end{tabular}
\]
Here the $\sqrt{An}$ arises from undoing the normalization and the
``2'' is the value of $L$.
Of course, this is \textit{rough}: we mean
``closer to $2$ than to $1$ or $3.$''

\subsection{Other Results for the Random Network Models}
There is substantial literature on the networks (MST, proximity graphs,
Delaunay triangulation) in the deterministic setting.
In the random case,
central limit theorems for total network length have been studied in
many models:
for the MST in \cite{KL96,MR1484795,MR1747451},
and for the Delaunay triangulation, Voronoi tessellation, relative
neighborhood and Gabriel graphs in
\cite{MR1241033,MR1272628,MR1878288}.
Large deviation estimates for total network length are given for the Gabriel
graph in \cite{MR1361756}, Section 11.4, and presumably could be
extended to other models.
Otherwise the literature for the random case is rather diffuse,
with different focuses for different networks.
For instance, work on MSTs has focused on
connections with critical continuum percolation \cite{MR1645627}.
For the relative neighborhood graph and the Gabriel graph,
\cite{MR937563} calculates $\bar{\Delta}$ and
\cite{MR2257270} shows that, in the finite model, in a certain range
the $\beta$-skeletons have
%
%
\begin{equation}
\Rmax\mbox{ grows as order }
\sqrt{\log n/\log\log n}
\label{loglog}
\end{equation}
and \cite{DGM09} shows the same order for maximum vertex degree in the
Gabriel graph.
As for the Delaunay triangulation, there has been surprisingly little
follow-up to the seminal analysis by Miles \cite{miles70} (various
maximal statistics are studied in \cite{MR1099499}), though the
closely related
Voronoi tessellation has been studied in more detail \cite{MR1295245}.

\subsection{Speculative Applications of Random Proximity Graphs}
\label{sec-specu}
Random proximity graphs seem an interesting object of study from many
viewpoints, in particular, as an attractive alternative to random
geometric graphs for modeling spatial networks that are connected by design.
It is remarkable that results such as (\ref{loglog}) are the only
nonelementary results about them that we can find in the literature.
As well as being natural models for road networks,
proximity graphs might be useful in modeling
communication networks suffering line of sight interference.

At a more mathematical level,
for questions such as spread-out percolation \cite{MR1202526} or
critical value of contact processes \cite{MR2298278},
random proximity graphs with small $A$ are an interesting alternative
to the usual lattice- or random graph-based models.
For instance, it is natural to conjecture that the critical value
$p^*_A$ for edge percolation on a random proximity graph with template $A$
satisfies
%
%
\begin{equation}
p^*_A \sim\pi^{-1} \operatorname{area}(A) \quad \mbox{as } \operatorname{area}(A)
\to0
\end{equation}
[the right side $= 1/\bar{\Delta}$ from (\ref{degree-scale})]
and that the critical value $\lambda^*_A$ for the contact process
has the same asymptotics.

\subsection{Analogies Between Deterministic and Random Networks}
\label{sec-analogs}
As mentioned earlier, we may make very loose analogies between
particular networks on random points and particular deterministic
networks in Figure~\ref{fig4}, based in part on exact equality of $\bar{\Delta}$
in the latter three cases:
\begin{eqnarray*}
\mbox{Relative n'hood graph} & \leftrightarrow& \mbox{punctured
lattice}, \\
\mbox{Gabriel graph} & \leftrightarrow& \mbox{square lattice}, \\
\mbox{Hammersley network} & \leftrightarrow& \mbox{diagonal lattice},
\\
\mbox{Delaunay triangulation} & \leftrightarrow& \mbox{triangular
lattice}.
\end{eqnarray*}

\subsection{Scale Invariant Continuum Networks}
Introducing the statistic $R$ can be viewed as one approach to
resolving the
``paradox'' from \cite{me116},
discussed in Section \ref{sec-Rs}, that the more natural statistic
$\Rave$
does not lead to realistic optimal networks in the $n \to\infty$ limit.
This particular approach was prompted by visualizing real-world road
networks---cf. discussion in Section \ref{sec-max23}.
Let us mention a mathematically more sophisticated alternative,
under study as a work in progress \cite{me-SIRSN}.
Instead of a discrete Poisson process of cities, we imagine a continuum
limit. That is,
for each finite set $(z_1,\ldots,z_k)$ of points in the plane, there
is a
random network $\mathcal{S}(z_1,\ldots,z_k)$ linking the points,
consistent as
more points are added.
Mathematically natural structural properties for
the distribution of such a process are as follows:

\begin{longlist}[(ii)]
\item[(i)] translation and rotation invariance,
\item[(ii)] scale invariance,
\end{longlist}
where the latter means that routes,
\textit{as point-sets in ${\mathbb{R}}^2$}, are invariant in
distribution under
Euclidean scaling.
This implies that the quantity $\rho(d)$ analogous to (\ref{R-def}),
assumed finite, is a constant, which we can call $R^\prime$. The analog
$L^\prime$ of $L$ is defined by
\[
\begin{tabular}{p{200pt}@{}}
the expected length of the network on $n$ uniform random points in the
area-$n$ square grows
$\sim L^\prime n$ as $n \to\infty$.
\end{tabular}
\]
In this setting we can study the optimal trade-off between $L^\prime$
and $R^\prime$, and the kind of
``paradoxical'' Figure \ref{fig5} network cannot arise because it violates
scale-invariance.

\section*{Acknowledgments}
Aldous's research supported by
NSF Grant DMS-0704159.
We thank three anonymous referees for helpful comments.


\begin{thebibliography}{10}

\bibitem{me118}
\textsc{Aldous, D. J.} (2008).
Spatial transportation networks with transfer costs:
Asymptotic optimality of hub and spoke models.
\textit{Math. Proc. Cambridge Philos. Soc.} \textbf{145} 471--487.
\MR{2442138}

\bibitem{me112}
\textsc{Aldous, D. J.} (2008).
Optimal spatial transportation networks where link-costs are
sublinear in link-capacity.
\textit{J. Stat. Mech.} \textbf{2008} P03006.

\bibitem{me-spatial-4}
\textsc{Aldous, D. J.} (2009).
Which connected spatial networks on random points have linear
route-lengths?
Available at \href{http://arxiv.org/abs/arXiv:0911.5296v1}{arXiv:0911.5296v1}.

\bibitem{me-spatial-5}
\textsc{Aldous, D. J.} (2009).
The shape theorem for route-lengths in connected spatial
networks on random
points.
Available at \href{http://arxiv.org/abs/arXiv:0911.5301v1}{arXiv:0911.5301v1}.

\bibitem{me-SIRSN}
\textsc{Aldous, D. J.} (2010).
Scale-invariant random spatial networks.
To appear.

\bibitem{me71}
\textsc{Aldous, D. J.} and \textsc{Diaconis, P.} (1995).
Hammersley's interacting particle process and longest increasing
subsequences.
\textit{Probab. Theory Related Fields} \textbf{103} 199--213.
\MR{1355056}

\bibitem{me116}
\textsc{Aldous, D. J.} and \textsc{Kendall, W. S.} (2008).
Short-length routes in low-cost networks via {P}oisson line
patterns.
\textit{Adv. in Appl. Probab.} \textbf{40} 1--21.
\MR{2411811}


\bibitem{me-spatial-3}
\textsc{Aldous, D. J., Bhamidi, S.} and \textsc{Lando, T.} (2010).
The stretch-length tradeoff in geometric networks: Worst-case and
average-case study.
To appear.

\bibitem{me57}
\textsc{Aldous, D. J.} and \textsc{Steele, J. M.} (1992).
Asymptotics for {E}uclidean minimal spanning trees on random points.
\textit{Probab. Theory Related Fields} \textbf{92} 247--258.
\MR{1161188}

\bibitem{me-spatial-2}
\textsc{Aldous, D. J.} and \textsc{Choi, A.} (2009).
A route-length\break efficiency statistic for road networks.
Unpublished manuscript. Available at
\href{http://www.stat.berkeley.edu/~aldous/Spatial/paper.pdf}{www.stat.berkeley.edu/}\break
\href{http://www.stat.berkeley.edu/~aldous/Spatial/paper.pdf}{\textasciitilde aldous/Spatial/\break paper.pdf}.

\bibitem{AB92}
\textsc{Avram, F.} and \textsc{Bertsimas, D.} (1992).
The minimum spanning tree constant in geometric probability
and under
the independent model: A unified approach.
\textit{Ann. Appl. Probab.} \textbf{2} 113--130.
\MR{1143395}

\bibitem{MR1241033}
\textsc{Avram, F.} and \textsc{Bertsimas, D.} (1993).
On central limit theorems in geometrical probability.
\textit{Ann. Appl. Probab.} \textbf{3} 1033--1046.
\MR{1241033}

\bibitem{BTZ00}
\textsc{Baccelli, F., Tchoumatchenko, K.} and \textsc{Zuyev, S.} (2000).
{Markov paths on the Poisson--Delaunay graph with applications to
routeing in mobile networks}.
\textit{Adv. in Appl. Probab.} \textbf{32} 1--18.
\MR{1765174}

\bibitem{barth06}
\textsc{Barth\'{e}lemy, M.} and \textsc{Flammini, A.} (2006).
Optimal traffic networks.
\textit{J. Stat. Mech. Theory Exp.}
\textbf{2006} L07002.

\bibitem{MR2298278}
\textsc{Berger, N., Borgs, C., Chayes, J. T.} and \textsc{Saberi, A.} (2005).
On the spread of viruses on the internet.
In \textit{Proceedings of the Sixteenth Annual ACM-SIAM
Symposium on
Discrete Algorithms} 301--310 (electronic). ACM, New York.
\MR{2298278}

\bibitem{MR1099499}
\textsc{Bern, M., Eppstein, D.} and \textsc{Yao, F.} (1991).
The expected extremes in a {D}elaunay triangulation.
\textit{Internat. J. Comput. Geom. Appl.} \textbf{1} 79--91.
\MR{1099499}

\bibitem{MR1645627}
\textsc{Bezuidenhout, C., Grimmett, G.} and \textsc{L{\"o}ffler, A.} (1998).
Percolation and minimal spanning trees.
\textit{J. Stat. Phys.} \textbf{92} 1--34.
\MR{1645627}

\bibitem{MR2257270}
\textsc{Bose, P., Devroye, L., Evans, W.} and \textsc{Kirkpatrick,
D.} (2006).
On the spanning ratio of {G}abriel graphs and {$\beta$}-skeletons.
\textit{SIAM J. Discrete Math.} \textbf{20} 412--427 (electronic).
\MR{2257270}

\bibitem{cortina}
\textsc{Cortina-Borja, M.} and \textsc{Robinson, T.} (2000).
Estimating the asymptotic constants of the total length of
{E}uclidean minimal spanning trees with power-weighted edges.
\textit{Statist. Probab. Lett.} \textbf{47} 125--128.
\MR{1747099}

\bibitem{MR937563}
\textsc{Devroye, L.} (1988).
The expected size of some graphs in computational geometry.
\textit{Comput. Math. Appl.} \textbf{15} 53--64.
\MR{0937563}

\bibitem{DGM09}
\textsc{Devroye, L., Gudmundsson, J.} and \textsc{Morin, P.} (2009).
On the expected maximum degree of Gabriel and Yao graphs.
\textit{Adv. in Appl. Probab.} \textbf{41} 1123--1140.
\MR{2663239}

\bibitem{gastner-newman}
\textsc{Gastner, M. T.} and \textsc{Newman, M. E. J.} (2006).
Shape and efficiency in spatial distribution networks.
\textit{J. Stat. Mech. Theory Exp.} \textbf{2006} P01015 (electronic).

\bibitem{MR1686530}
\textsc{George, P.-L.} and \textsc{Borouchaki, H.} (1998).
\textit{Delaunay Triangulation and Meshing}.
Editions Herm\`es, Paris.
\MR{1686530}

\bibitem{graham-hell}
\textsc{Graham, R. L.} and \textsc{Hell, P.} (1985).
On the history of the minimum spanning tree problem.
\textit{IEEE Ann. History Comput.} \textbf{07} 43--57.
\MR{0783327}

\bibitem{MR1272628}
\textsc{Heinrich, L.} (1994).
Normal approximation for some mean-value estimates of absolutely
regular tessellations.
\textit{Math. Methods Statist.} \textbf{3} 1--24.
\MR{1272628}

\bibitem{Poisson-Matching}
\textsc{Holroyd, A. E., Pemantle, R., Peres, Y.} and \textsc
{Schramm,~O.} (2009).
Poisson matching.
\textit{Ann. Inst. H. Poincar\'e Probab. Statist.} \textbf{45} 266--287.
\MR{2500239}

\bibitem{MR1961286}
\textsc{Holroyd, A. E.} and \textsc{Peres, Y.} (2003).
Trees and matchings from point processes.
\textit{Electron. Comm. Probab.} \textbf{8} 17--27 (electronic).
\MR{1961286}

\bibitem{jaromczyk}
\textsc{Jaromczyk, J. W.} and \textsc{Toussaint, G. T.} (1992).
Relative neighborhood graphs and their relatives.
\textit{Proc. IEEE} \textbf{80} 1502--1517.

\bibitem{KL96}
\textsc{Kesten, H.} and \textsc{Lee, S.} (1996).
The central limit theorem for weighted minimal spanning trees on
random points.
\textit{Ann. Appl. Probab.} \textbf{6} 495--527.
\MR{1398055}

\bibitem{kirk-radke}
\textsc{Kirkpatrick, D. G.} and \textsc{Radke, J. D.} (1985).
A framework for computational morphology.
In \textit{Computational Geometry} (G.~T.~Toussaint, ed.)
217--248. Elsevier, Amsterdam.

\bibitem{MR1484795}
\textsc{Lee, S.} (1997).
The central limit theorem for {E}uclidean minimal spanning trees.
{I}.
\textit{Ann. Appl. Probab.} \textbf{7} 996--1020.
\MR{1484795}

\bibitem{MR1747451}
\textsc{Lee, S.} (1999).
The central limit theorem for {E}uclidean minimal spanning trees.
{II}.
\textit{Adv. in Appl. Probab.} \textbf{31} 969--984.
\MR{1747451}

\bibitem{MR0660693}
\textsc{Lloyd, E. L.} (1977).
On triangulations of a set of points in the plane.
In \textit{18th Annual Symposium on Foundations of Computer Science
(Providence, RI, 1977)} 228--240. IEEE Comput. Soc., Long Beach,
CA.
\MR{0660693}

\bibitem{matula}
\textsc{Matula, D. W.} and \textsc{Sokal, R. R.} (1980).
Properties of {G}abriel graphs relevant to geographical
variation research and the clustering of points in the plane.
\textit{Geog. Anal.} \textbf{12} 205--222.

\bibitem{miles70}
\textsc{Miles, R. E.} (1970).
On the homogeneous planar {P}oisson point process.
\textit{Math. Biosci.} \textbf{6} 85--127.
\MR{0279853}

\bibitem{MR1295245}
\textsc{M{\o}ller, J.} (1994).
\textit{Lectures on Random {V}orono\u\i\ Tessellations}.
\textit{Lecture Notes in Statistics} \textbf{87}.
Springer, New York.
\MR{1295245}

\bibitem{MR2289615}
\textsc{Narasimhan, G.} and \textsc{Smid, M.} (2007).
\textit{Geometric Spanner Networks}.
Cambridge Univ. Press, Cambridge.
\MR{2289615}

\bibitem{kumar02}
\textsc{Narayanaswamy, S., Kawadia, V., Sreenivas, R. S.} and \textsc
{Kumar, P. R.} (2002).
Power control in ad-hoc networks: Theory, architecture,
algorithm and
implementation of the {COMPOW} protocol.
In \textit{Proc. European Wireless Conference, Florence, Italy}.

\bibitem{newman-survey}
\textsc{Newman, M. E. J.} (2003).
The structure and function of complex networks.
\textit{SIAM Rev.} \textbf{45} 167--256.
\MR{2010377}

\bibitem{penrose-RGG}
\textsc{Penrose, M. D.} (2003).
\textit{Random Geometric Graphs}.
Oxford Univ. Press, Oxford.
\MR{1986198}

\bibitem{MR1202526}
\textsc{Penrose, M. D.} (1993).
On the spread-out limit for bond and continuum percolation.
\textit{Ann. Appl. Probab.} \textbf{3} 253--276.
\MR{1202526}

\bibitem{MR1878288}
\textsc{Penrose, M. D.} and \textsc{Yukich, J. E.} (2001).
Central limit theorems for some graphs in computational geometry.
\textit{Ann. Appl. Probab.} \textbf{11} 1005--1041.
\MR{1878288}

\bibitem{MR1952000}
\textsc{Penrose, M. D.} and \textsc{Yukich, J. E.} (2003).
Weak laws of large numbers in geometric probability.
\textit{Ann. Appl. Probab.} \textbf{13} 277--303.
\MR{1952000}\

\bibitem{steele97}
\textsc{Steele, J. M.} (1997).
\textit{Probability Theory and Combinatorial Optimization}.
\textit{CBMS-NSF Regional Conference Series in Applied Math} \textbf{69}.
SIAM, Philadelphia, PA.
\MR{1422018}

\bibitem{MR895588}
\textsc{Stoyan, D., Kendall, W. S.} and \textsc{Mecke, J.} (1995).
\textit{Stochastic Geometry and Its Applications}, 2nd ed.
Wiley, Chichester.
\MR{0895588}

\bibitem{MR1361756}
\textsc{Talagrand, M.} (1995).
Concentration of measure and isoperimetric inequalities in product
spaces.
\textit{Inst. Hautes \'Etudes Sci. Publ. Math.} \textbf{81} 73--205.
\MR{1361756}

\bibitem{yukich-book}
\textsc{Yukich, J. E.} (1998).
\textit{Probability Theory of Classical {E}uclidean Optimization
Problems}.
\textit{Lecture Notes in Math.} \textbf{1675}. Springer, Berlin.
\vspace*{-2pt}
\MR{1632875}
\end{thebibliography}
\end{document}